\documentstyle[12pt]{article}

\catcode `\@=11
\@addtoreset{equation}{section}
\def\theequation{\arabic{section}.\arabic{equation}}
\catcode `\@=12

\newcommand{\be}{\begin{equation}}
\newcommand{\en}{\end{equation}}
\newcommand{\bea}{\begin{eqnarray}}
\newcommand{\ena}{\end{eqnarray}}
\newcommand{\beano}{\begin{eqnarray*}}
\newcommand{\enano}{\end{eqnarray*}}
\newcommand{\bee}{\begin{enumerate}}
\newcommand{\ene}{\end{enumerate}}
\newcommand{\Bbb}{\rm\bf}
\newcommand{\R}{R \!\!\!\! R}
\newcommand{\N}{N \!\!\!\!\! N}
\newcommand{\Z}{Z \!\!\!\!\! Z}
\newcommand{\Hil}{{\cal H}}
\newcommand{\I}{{\cal I}}

\newcommand{\Id}{1\!\!1}
\newcommand{\F}{{\cal F}}
\newcommand{\Log}{\mbox{Log}}

\begin{document}

\thispagestyle{empty}
 
\vspace*{1cm}

\begin{center}
{\Large \bf Some Results About Frames}   \vspace{2cm}\\

{\large F. Bagarello}
\vspace{3mm}\\
  Dipartimento di Matematica ed Applicazioni, 
Fac. Ingegneria, Universit\`a di Palermo, I-90128  Palermo, Italy\\
e-mail: Bagarello@ipamat.math.unipa.it
\vspace{2mm}\\
\end{center}

\vspace*{2cm}

\begin{abstract}
\noindent 
In this paper we discuss some topics related to the general theory of frames. In particular
we focus our attention to the existence of different 'reconstruction formulas' for a given
vector of a certain Hilbert space and to some refinement of the perturbative approach for
the computation of the dual frame.  \end{abstract}

\vfill

\newpage

\section{Introduction}
Whenever we deal with a (separable) Hilbert space $\Hil$ the first problem we face with is
the way in which an arbitrary element $f \in \Hil$ can be conveniently expressed. As we
know the usual choice is to expand $f$ in terms of an orthonormal basis $\{e_n\}$ of
$\Hil$: in this way the expansion is particularly simple, $f=\sum_n<e_n,f>e_n$, and the
Parseval equality is satisfied, $\sum_n|<e_n,f>|^2=\|f\|$.
It is also well known that in any Hilbert space we can build up different orthonormal
bases, for instance using unitary transformations starting with one of such basis.

Sometimes, however, the conditions of our (physical) problem force us to consider a set of
vectors $\{\Phi_n\}$ which is no longer orthonormal but is still a basis of $\Hil$. We have
a typical example of this situation when the set $\{\Phi_n\}$ forms what is known as
Riesz basis, see \cite{cdf} and references therein. One possible characterization of a
Riesz basis is the following:

The set $\{\Phi_n\}$ is a Riesz basis if and only if

- all the $\Phi_n$ are independent, that is no $\Phi_{n_o}$ lies within the closure of
the finite linear span of the other $\Phi_n$, and

- $\exists \, A>0, \: B<\infty$ so that, for any $f \in \Hil$,
\be
A\|f\|^2 \leq \sum_n|<\Phi_n,f>|^2 \leq B\|f\|^2.
\label{11}
\en
This last property implies that the vectors of the set $\{\Phi_n\}$ generates the
whole Hilbert space, while the first condition says that these vectors are linearly
independent: therefore they form a {\em typical} basis of $\Hil$. 

We know many examples of sets of vectors which are not a Riesz basis but still have a
relevant role in the description of vectors of $\Hil$. For instance, if we consider any
overcomplete set of coherent states, \cite{kla}, this set satisfies a relation similar to
the one in (\ref{11}), but the vectors are not linearly independent. Equation (\ref{11}) is
also satisfied by some sets of wavelets, \cite{dau1}. Sets of vectors of this kind are
known as frames. In other words, we can say that a frame is a set of generators of
$\Hil$, but, since the vectors are not independent, the way in which a vector $f\in
\Hil$ can be expanded by means of these vectors is not unique.

Many properties of the frames have been discussed in the literature. In particular the
interest of the researchers has been focused mainly on the reconstruction formula: given
a frame  $\{\varphi_n\}$ in which way the sequence $\{<\varphi_n,f>\}$ is linked to the
function $f$? What does the Parseval equality become for frames? It is possible to
operate in such a way to optimize the numerical procedure?

All these questions have already an answer in the literature. In this paper we will show
that many other relevant questions are still to be raised and can be solved: it is
possible to find different reconstruction formulas, maybe more convenient than the usual
one? It is possible to build up a perturbation scheme which works in conditions in
which the usual one, \cite{dau1}, has difficulties?

The paper is organized as follows:

in Section 2 we introduce a bit of notation and, in order to keep the paper self
contained, we quickly recall the main known results about frames.

In Section 3 we show how, starting from a given frame, it is always possible to define a
tight frame with frame bound equal to one, as well as many other frames which give
different reconstruction formulas. We will discuss many examples of this procedure and we
will also show how the usual perturbation scheme works in this case.

In Section 4 we discuss a different perturbation approach which may be useful whenever the
frame constants are not close to each other.

In Section 5 we discuss some analogies between frames and orthonormal bases.

\section{Notation and Known Results}

Most of the results we are going to summarize in this Section can be found
in \cite{dau1} and \cite{dau2}. We also like to cite \cite{kol} which is an useful review
on wavelets with a long Section on frames and where it is contained an excellent
bibliography.

Let $\Hil$ be an Hilbert
space and $J$ a given set of indexes. Let also $\I\equiv \{\varphi_n(x)\}$, $n\in J$, be a
set of vectors of $\Hil$. We say that $\I$ is an (A,B)-{\em frame} of $\Hil$ if there
exist two positive constants, called frame bounds, $A>0$ and $B<\infty$, such that the
inequalities\be A\|f\|^2 \leq \sum_n|<\varphi_n,f>|^2 \leq B\|f\|^2
\label{21}
\en
hold for any $f\in \Hil$.

To any such set $\I$ can be associated a bounded operator $F:\Hil \rightarrow {\Bbb
l}^2(J)=\{\{c_n\}_{n\in J}\, : \: \sum_{n\in J}|c_n|^2\}<\infty$  defined by the formula
\be
\forall f\in \Hil \hspace{1cm}  (Ff)_j = <\varphi_j, f>.
\label{22}
\en
We will omit the dependence on $\I$ of $F$ whenever this does not cause confusion.
Due to equation (\ref{21}) we see that $\|F\| \leq \sqrt{B}$. The
adjoint of the operator $F$, $F^*$, which maps ${\Bbb l}^2(J)$ into $\Hil$, is such that
 \be
\forall \{c\} \in {\Bbb l}^2(J) \hspace{1cm}  F^*c \equiv \sum_{i\in J}c_i \varphi_i
\label{23}
\en
Of course, also $F^*$ depends on the family $\I$ we are considering. Again, we will omit
this dependence whenever it will not cause confusion.

By means of these operators condition (\ref{21}) can be rewritten in the following
equivalent way: $\I$ is an (A,B)-{\em frame} of $\Hil$ if there exist
two positive constants, $A>0$ and $B<\infty$, such that the inequalities
\be
A\Id \leq F^* F \leq B \Id
\label{24}
\en
hold in the sense of the operators, \cite{rs}. We have used $\Id$ to identify the
identity operator in $B(\Hil)$.

Condition (\ref{24}) implies that the operator $F^* F$, which maps $\Hil$ into itself, can
be inverted and that its inverse, $(F^* F)^{-1}$ is still bounded in $\Hil$. In other
terms, we have that both $F^*F$ and $(F^*F)^{-1}$ belong to $B(\Hil)$.

Following the literature, see \cite{dau1} for instance, one define the {\em dual frame} of
$\I$, $\tilde \I$,  the set of vectors $\tilde \varphi_i$ defined by
\be
\tilde \varphi_i \equiv (F^* F)^{-1} \varphi_i \hspace{2cm} \forall i\in J.
\label{25}
\en
We have called  $\tilde \I$ {\em the dual frame} since, as a matter of fact, it is really a
frame, and, in particular, it is a $(\frac{1}{B},\frac{1}{A})$-frame, see \cite{dau1}.
Defining now a new operator between $\Hil$ and ${\Bbb l}^2(J)$ as $\tilde F \equiv F
(F^*F)^{-1}$, it is easy to prove that $\tilde F$ is such that $(\tilde F f)_i =<\tilde
\varphi_i, f>$, for all $f\in \Hil$. Moreover, the following relations hold: $\tilde F^*
F=F^*\tilde F =\Id$. By means of these
 equalities it is proved that any vector of the Hilbert space can be expanded as 
linear combinations of the vectors of the set $\I$ or of the set $\tilde \I$. We have the
following reconstruction formulas: 
\be
f=\sum_{i\in J}<\varphi_i,f> \tilde \varphi_i = \sum_{i\in J}<\tilde \varphi_i,f> \varphi_i
\label{26}
\en
for all $f\in \Hil$. In \cite{dau1} it is also discussed that, since
$\varphi_i= \tilde{\tilde \varphi_i}$ for all $i\in J$, then the dual frame of the set
$\tilde \I$ is nothing but the set $\I$ itself.

As it is clear from equation (\ref{26}), a crucial role in the reconstruction procedure is
the knowledge of the set $\tilde \I$. In order to obtain the explicit expression for $\tilde
\varphi_i$, we first have to know how the operator $(F^*F)^{-1}$ acts on the vectors of
$\Hil$. This is, in general, a difficult problem to solve. Only in a single situation we can
give an easy answer, namely when our frame is {\em tight}. This means that the frame bounds
$A$ and $B$ coincide, $A=B$, so that equation (\ref{24})
reduces to  
\be
 F^* F = A \Id,
\label{27}
\en
which implies also that $(F^*F)^{-1}=\frac{1}{A} \Id$. Therefore $\tilde \varphi_i=
\frac{1}{A} \varphi_i$, for all $i\in J$. In this case the reconstruction formulas above
coincide and they look like
$$
f=\frac{1}{A} \sum_{i\in J}<\varphi_i,f> \varphi_i,
$$
for all $f\in \Hil$. In particular, moreover, if $A=1$ and if all the vectors $\varphi_i$
are normalized, it follows that this frame forms an orthonormal basis of the Hilbert
space, see \cite{dau1}. Of course, also the vice-versa holds true: if $\I$ is an orthonormal
set in $\Hil$, then $\I$ is a $(1,1)$-frame of normalized vectors. This is an obvious
consequence of the Parceval equality $\sum_n|<\varphi_n,f>|^2 =\|f\|^2$, which holds for all
$f\in \Hil$.

How far can we go if the frame we are dealing with is not tight? In \cite{dau1} and in the
references therein, it is proposed a perturbative approach to this problem, which works well
whenever $A\approx B$. We will discuss in Section 4 a different approximation approach which
can be used even when the frame bounds $A$ and $B$ are not very close. In the final
part of this Section we summarize the standard procedure. We start defining a bounded
self-adjoint operator $R$ by \be
R=\Id -\frac{2}{A+B}F^*F.
\label{28}
\en
Using inequality (\ref{24}) we find that 
\be
-\frac{B-A}{B+A} \leq R \leq \frac{B-A}{B+A},
\label{29}
\en
which implies that
\be
\|R\| \leq \frac{B-A}{B+A}.
\label{210}
\en
Inverting definition (\ref{28}) we obtain $F^*F=\frac{A+B}{2}(\Id-R)$, which also implies
that  $(F^*F)^{-1}=\frac{2}{A+B}(\Id-R)^{-1}$. Since certainly $\|R\|<1$ we can expand the
operator $(\Id-R)^{-1}$ into a norm-convergent series, $\sum_{k=0}^{\infty}R^k$. Therefore
\be
\tilde \varphi_j = \frac{2}{A+B}\sum_{k=0}^{\infty}R^k \varphi_j.
\label{211}
\en
The $N-th$ approximation is obtained considering only the first $N-th$ contributions in
the above infinite sum: 
\be
\tilde \varphi_j^N = \frac{2}{A+B}\sum_{k=0}^{N}R^k \varphi_j.
\label{212}
\en
In \cite{dau1} it is shown that, for all $f\in \Hil$, the following estimate holds:
\be
\left\|f-\sum_{j\in J}<\varphi_j,f>\tilde \varphi_j^N\right\| \leq \|R\|^{N+1}\|f\|,
\label{213}
\en
which says that the above quantity converges to zero when $N\rightarrow \infty$. Of
course, due to (\ref{210}), this convergence is fast when $A\approx B$, while it is rather
slow if $B-A$ is big, so that we cannot approximate $\tilde \varphi_j$ keeping only few
contributions in the expansion (\ref{211}). We need to consider $N$ big enough: this is
the only way in which the error on the right hand side of (\ref{213}) can be made as small
as we like. A possibility to overcome this problem will be proposed in Section 4.

\section{More Reconstruction Formulas}

We have seen in the previous Section that, given an $(A,B)$-frame $\I\equiv
\{\varphi_n\}$, it is possible to associate to it an unique dual
$(\frac{1}{B},\frac{1}{A})$-frame $\tilde \I\equiv \{\tilde \varphi_n\}$ and to use this
two frames {\em together} to expand any vector of the Hilbert space. As it is well known,
this is only one of infinite equivalent possibilities: in fact, since the sets $\I$ and
$\tilde \I$ are not made up of independent vectors, in general, we cannot invoke any
unicity of the coefficients of the expansion of a given vector. In other words: given a
vector $f\in \Hil$ and a frame $\I$, in general there exist infinitely many way to expand
$f$ in terms of $\varphi_j$. Of course the possibilities given by (\ref{26}) are
particularly simple and elegant, so that they are considered of a particular interest.
Moreover, these expansions solve the problem raised in the Introduction about the relation
between the sequence $\{<\varphi_n,f>\}$ and the vector $f$, $\{\varphi_n\}$ being a frame.

The problem we want to deal with in this Section is the following: does a given frame
suggests other resolutions of the identity or, which is the same, other reconstruction
formulas?

The answer to this question is yes, and we will show how these different reconstruction
formulas can be obtained and which techniques must be used. In particular we will show that
the procedure discussed in the previous Section is only a particular case of this more
general approach.

Our first step consists in defining the operator
\be
\F_1\equiv F^*F.
\label{31}
\en
The norm of this operator is bounded from above and from below, $A\leq \|\F_1\| \leq B$, it
is positive, $\F_1 \geq 0$, self adjoint, $\F_1 = \F_1^*$, and its action on a given vector
of $\Hil$ is given by $\F_1f=\sum_{i\in J}<\varphi_i,f>\varphi_i$.

Let $E_\lambda$ be the family of spectral operators of $\F_1$. We can write, making
use of the spectral theorem,
\be
<\F_1\Phi,\Psi> = \int_A^B\lambda \: d<E_\lambda\Phi, \Psi>, \hspace{2cm} \forall \Phi, \Psi
\in \Hil.
\label{32}
\en
Due to the fact that $0<A\leq B<\infty$, we can define arbitrary powers, positive and
negative, of the operator $\F_1$:
\be
<\F_\alpha \Phi,\Psi> \equiv <(F^*F)^\alpha\Phi,\Psi>  = \int_A^B\lambda^\alpha
\,d<E_\lambda \Phi, \Psi>,  
\label{33}
\en
for all $\Phi, \Psi \in \Hil$, and $\forall \alpha \in \R$.

Following \cite{ric} we see that $\F_\alpha$, which of course still maps $\Hil$ into itself,
is self adjoint and bounded for any $\alpha \in \R$. Moreover we have:
\be
(\F_\alpha)^\beta =\F_{\alpha\beta}, \hspace{2cm} \F_\alpha \F_\beta = \F_{\alpha+\beta}.
\label{34}
\en
Using the spectral theorem it is an easy exercise to compute the bounds for the operator
$\F_\alpha$. We find that
\bea
&&A^\gamma \Id\leq \F_\gamma \leq B^\gamma \Id \hspace{1cm} \forall \gamma \geq 0 
\label{35} \\
&& B^\gamma \Id \leq \F_\gamma \leq A^\gamma \Id \hspace{1cm} \forall \gamma < 0.
\label{36}
\ena
Given an arbitrary real number $\alpha$ let us define the following vectors:
\be
\varphi_i^{(\alpha)} \equiv \F_\alpha \varphi_i \hspace{2cm} \forall i\in J,
\label{37}
\en
and let us call $\I^{(\alpha)}$ the set of these vectors. We can prove that all these sets
are frames in $\Hil$. We have indeed the following

\noindent
	{\bf Proposition 1}.

 Each $\I^{(\alpha)}$ is a frame. In particular $\I^{(\alpha)}$ is an
$(A^{2\alpha+1},B^{2\alpha+1})$-frame if $\alpha>-\frac{1}{2}$, is a
$(1,1)$-frame if $\alpha=-\frac{1}{2}$, and is a $(B^{2\alpha+1},A^{2\alpha+1})$-frame if
$\alpha<-\frac{1}{2}$.

\vspace{5mm}

\noindent
	{\underline {Proof}}

This result easily follows from inequalities (\ref{35}) and (\ref{36}) and by the following
equalities:
\beano
&&\sum_{i\in J}|<\varphi_i^{(\alpha)},f>|^2= \sum_{i\in J}|<\F_\alpha \varphi_i,f>|^2= 
\sum_{i\in J}|<\varphi_i,\F_\alpha f>|^2 = \\
&&= \|F(\F_\alpha f)\|^2=<\F_\alpha \F_1 \F_\alpha f,f> = <\F_{2\alpha+1} f,f>,
\enano
which hold for all vectors $f\in \Hil$.\hfill $\Box$

\vspace{3mm}

\noindent
	{\bf Remark}.-- The above procedure teaches, in particular, how to get a tight frame with
frame bound equal to 1 starting from a generic frame. Of course this does not imply that the
vectors $\varphi_i^{(-1/2)}$ form an orthonormal basis since normalization of these vectors
is not ensured.

\vspace{3mm}

We now define an operator $\F^{(\alpha)}$ 'associated' to the frame $\I^{(\alpha)}$, operator
which plays the same role of $\tilde F$ in the standard case. 

Let $\F^{(\alpha)}$ be an operator from $\Hil$ into ${\Bbb l}^2(J)$ 
defined by the formula
\be
\forall f\in \Hil \hspace{1cm}  (\F^{(\alpha)}f)_j = <\varphi_j^{(\alpha)}, f>.
\label{38}
\en
Since we have also $<\varphi_j^{(\alpha)}, f>= <\F_\alpha \varphi_j, f> = <\varphi_j
,\F_\alpha f> = (F(\F_\alpha f))_j$, it follows that
\be
\F^{(\alpha)} = F \F_\alpha.
\label{39}
\en
Using now the equalities in (\ref{34}) we see that, if we define the operator
\be
\tilde \F^{(\alpha)} \equiv \F_{-1-\alpha} F^*,
\label{310}
\en
then
\be
\tilde \F^{(\alpha)} \F^{(\alpha)} = (\F^{(\alpha)})^* (\tilde \F^{(\alpha)})^* =\Id,
\label{311}
\en
for all $\alpha \in \R$. These equalities can now be used to built up different
reconstruction formulas, one for each value of $\alpha$. First of all we observe that
definition (\ref{38}) implies that \be
(\F^{(\alpha)})^* c = \sum_{i\in J} c_i \varphi_i^{(\alpha)},
\label{312}
\en
for all $c \in {\Bbb l}^2(J)$. From (\ref{310}) we have that, given $f\in \Hil$,
\beano
&&(\tilde \F^{(\alpha)})^*f=F(\F_{-1-\alpha}f) =
\left\{<\varphi_j,\F_{-1-\alpha}f>\right\}_{j\in J} = \\
&& = \left\{<\F_{-1-\alpha} \varphi_j,f>\right\}_{j\in J} = \left\{<
\varphi_j^{(-1-\alpha)},f>\right\}_{j\in J}.
\enano
Therefore, given $c \in {\Bbb l}^2(J)$, we have
\be
\tilde \F^{(\alpha)}c = \sum_{i\in J} c_i \varphi_i^{(-1-\alpha)}.
\label{314}
\en
All these results can be collected to show that, for any $f\in \Hil$, we can write the
following expansions:
\be
f = \sum_{i\in J} <\varphi_i^{(\alpha)},f> \varphi_i^{(-1-\alpha)} = 
\sum_{i\in J} <\varphi_i^{(-1-\alpha)},f> \varphi_i^{(\alpha)},
\label{315}
\en
equalities which hold for all choices of the real $\alpha$ (Of course the second
expansion is redundant since it can be obtained from the first one putting $\alpha
\rightarrow -1-\alpha$).

The proof is very easy: using formulas (\ref{38}), (\ref{311}) and (\ref{314}) we see that
for any $f\in \Hil$, $f= \tilde \F^{(\alpha)} \F^{(\alpha)} f = \tilde \F^{(\alpha)}
(\F^{(\alpha)} f) = \sum_{i\in J} <\varphi_i^{(\alpha)},f> \varphi_i^{(-1-\alpha)}$. The
other equality is proved in the same way, by considering the second equality in (\ref{311}).

\vspace{8mm}

Let us now summarize our procedure:

(a) the starting point is a given frame $\I=\{\varphi_j\}$ and the related frame operators
$F$, $F^*$ and $\F_1=F^*F$;

(b) we construct the operator $\F_\alpha$ for a given real $\alpha$ as in (\ref{33});

(c) we define a new frame whose vectors are $\varphi_i^{(\alpha)}= \F_\alpha \varphi_i$;

(d) we also define its dual frame: $\varphi_i^{(-1-\alpha)}= \F_{-1-\alpha} \varphi_i$;

(e) we can now finally expand any $f\in \Hil$ as in (\ref{315}).

\vspace{4mm}

\noindent
	{\bf Remarks}.-- 

(a) It is interesting to notice that repeating twice the operation of taking
the dual of a given frame we come back to the original frame. This result generalizes the
analogous result for the standard procedure, $\varphi_j=\tilde{\tilde \varphi_j}$,
\cite{dau1}. More in details: if we take a frame
$\{\varphi_j\}$ and then we consider the new frame $\{\varphi_j^{(\alpha)}\}$ and its frame
operator $\F^{(\alpha)}$, it follows that the vectors constructed as in (\ref{25}),
$(\F^{(\alpha)*}  \F^{(\alpha)})^{-1} \varphi_j^{(\alpha)}$, coincide with
$\varphi_j^{(-1-\alpha)}$. In this sense we say that $\{\varphi_j^{(-1-\alpha)}\}$ is the
dual set of $\{\varphi_j^{(\alpha)}\}$. Moreover, if we consider the "bi-dual" of
$\{\varphi_j^{(\alpha)}\}$, then we obtain $\{\varphi_j^{(\alpha)}\}$ itself.

(b) The method proposed here looks as a concrete version of a more abstract approach
discussed in \cite{ant} where the attention was focused mainly on coherent states.

(c) We observe that the approach we have followed here can be seen as a special case of a
much more general one which is the following: given a frame $\{\varphi_j\}$ and its frame
operators $F$ and $\F_1$, considering the spectral decomposition (\ref{32}), we can define
new operators $F_f\equiv \int_A^B f(\lambda) \, dE_\lambda$, to be intended weakly, where
$f$ is any real bounded function with inverse bounded in $[A,B]$. In this way we could
define a new frame $\{\varphi_j^{(f)}\}=\{F_f \varphi_j\}$ and work on analogous
computations as the ones we have just discussed here, where we have always chosen
$f(\lambda)=\lambda^\alpha$, for some real $\alpha$.

\vspace{3mm}

Before discussing some examples of this procedure let us now consider some limiting case. We
begin with the simplest one, $\alpha =0$.

In this case we have $\varphi_i^{(\alpha)} = \F_0\varphi_i=\varphi_i$  and 
$\varphi_i^{(-1-\alpha)} = \varphi_i^{(-1)}= \F_{-1}\varphi_i=\tilde \varphi_i$, for all
$i\in J$. In this case the reconstruction formulas (\ref{315}) collapse to the standard
ones, (\ref{26}). 

Moreover the operators defined in this Section become: $\F^{(\alpha)}=F$, $\tilde
\F^{(\alpha)}=(F^*F)^{-1}F^*=\tilde F^*$, and so on. The equalities in (\ref{311}) coincide
now with the analogous ones $\tilde F^* F=F^*\tilde F =\Id$, discussed in the previous
Section in the standard setting.

\vspace{3mm}

The situation is completely specular for $\alpha =-1$. In this case each
$\varphi_i^{(\alpha)}$ coincides with $\tilde \varphi_i$ while $\varphi_i^{(-1-\alpha)}$ is
nothing but $\varphi_i$ itself. Again the operators $\F^{(\alpha)}$ and $\tilde
\F^{(\alpha)}$ coincide with the ones widely discussed in the literature and in Section 2.

\vspace{3mm}

Another limiting case is suggested by the Proposition 1 which assigns a particular
importance to the choice $\alpha =-\frac{1}{2}$. In this case we get a frame
$\I^{(-1/2)}=\{\varphi_i^{(-1/2)}\}$ which is tight with frame bound equal to 1. In fact,
if $\alpha=-\frac{1}{2}$, then $-1-\alpha=-\frac{1}{2}$. For any vector $f\in \Hil$ the
following equalities can be written
$$
f = \sum_{i\in J} <\varphi_i^{(-1/2)},f>
\varphi_i^{(-1/2)}, $$
and
$$
 \sum_{i\in J} |<\varphi_i^{(-1/2)},f>|^2 =\|f\|^2.
$$
For what concerns the expressions of the operators we have, for instance,
$\F^{(-1/2)}=F(F^*F)^{(-1/2)}$ and $\tilde \F^{(-1/2)} = (F^*F)^{(-1/2)} F^*=
(\F^{(-1/2)})^*$.

\vspace{3mm}

This last result is general: it is easy to show that for all $\alpha \in \R$ we have
\be
\F^{(\alpha)} = (\tilde \F^{(-1-\alpha)})^*.
\label{316}
\en

Let us now consider some examples. Of course, the situation is much simpler for finite
dimensional Hilbert spaces and it is just this situation that we start to consider. 

\vspace{4mm}

\noindent
{\bf Example 1} 

Let us take $\Hil = \R^2$, and let us consider the three vectors $\varphi_1=(1,0)$,
$\varphi_2=(0,1)$ and $\varphi_3=\frac{1}{\sqrt{2}}(1,1)$. It is easy to show that these
vectors satisfy the following relation
$$
\|f\|^2 \leq \sum_{i=1}^3|<\varphi_i,f>|^2 \leq 2\|f\|^2,
$$
for all $f\in \Hil$, so that $\I$ is a $(1,2)$-frame (the values of the bounds may be
optimized, but we don't care about this point here). In the perspective of showing in which
way our procedure works, it is convenient to start with frames which are not tight. This
is because we just want to show how a (1,1)-tight frame can be obtained starting with an
arbitrary other frame. It is obvious that $\I$  cannot be a Riesz basis, since its three
vectors are necessarily linearly dependent. 

For any $f=(f_x,f_y)$  we get
$\F_1f=\sum_{i=1}^3<\varphi_i,f>\varphi_i = \frac{1}{2}(3f_x+f_y,f_x+3f_y)$. Therefore, the
operator $\F_1$ can be identified with the $2\times 2$ matrix $$
\F_1\equiv  \frac{1}{2}\left(
\begin{array}{cc}
3 & 1  \\ 
1 & 3  \\ 
\end{array}
\right).
$$
Using a standard procedure we can write $\F_1$ in the more convenient form obtained by its
spectral decomposition:
\be
\F_1 = E_1+2E_2,
\label{317}
\en
where $\lambda_1=1$ and $\lambda_2=2$ are the eigenvalues of $\F_1$, the $E_k$'s are the
projection operators defined by $E_kg=<\eta_k,g>\eta_k$, $k=1,2$, $\forall g \in \Hil$, and
$\eta_1=\frac{1}{\sqrt{2}}(1,-1)$ and $\eta_2=\frac{1}{\sqrt{2}}(1,1)$ are the eigenstates
of $\F_1$ which correspond to $\lambda_1$ and $\lambda_2$.

We are now ready to build up different pairs of dual frames, each for any given value of
$\alpha$. Let us start with $\alpha=-1$. 

We have $\F_{-1} = E_1+\frac{1}{2}E_2$ and the set $\I^{(-1)}=\{\F_{-1}\varphi_i\}$ is
formed by the following vectors: $\varphi_1^{(-1)}=\frac{1}{4}(3,-1)$,
$\varphi_2^{(-1)}=\frac{1}{4}(-1,3)$ and $\varphi_3^{(-1)}=\frac{1}{2\sqrt{2}}(1,1)$. As we
have already discussed, $\I^{(-1)}$ is the dual frame of the original one $\I$. In the
standard language we would have written $\varphi_i^{(-1)} = \tilde \varphi_i$. It is
straightforward to verify that, given $f\in \Hil$, we can write $f=\sum_{i=1}^3<\varphi_i,f>
\varphi_i^{(-1)} = \sum_{i=1}^3<\varphi_i^{(-1)},f> \varphi_i$. Moreover we see that 
$\frac{1}{2}\|f\|^2 \leq \sum_{i=1}^3|<\varphi_i^{(-1)},f>|^2 \leq \|f\|^2$, estimate which
does agree with the statement of Proposition 1.

Of course the same results are obtained, mutatis mutandis, if we take $\alpha=0$. 

Let us now choose $\alpha=-\frac{1}{2}$. 

The operator $\F_{-\frac{1}{2}}$ can be simply deduced by (\ref{317}): $\F_{-\frac{1}{2}}= 
E_1+\frac{1}{\sqrt{2}}E_2$, so that we get the following vectors for the frame:
$\varphi_1^{(-1/2)}=\frac{1}{2}(1+\frac{1}{\sqrt{2}},-1+\frac{1}{\sqrt{2}})$,
$\varphi_2^{(-1/2)}=\frac{1}{2}(-1+\frac{1}{\sqrt{2}},1+\frac{1}{\sqrt{2}})$ and
$\varphi_3^{(-1/2)}=\frac{1}{2}(1,1)$. It is a simple calculation to verify that these
vectors form a $(1,1)$-frame. They are not, however, an orthonormal basis of $\Hil$ since
none of these vectors is normalized. This, of course, is what we must have
since the three vectors are necessarily linearly dependent in $\Hil$.

We conclude this example considering the case $\alpha=-\frac{1}{3}$.

The expression of the operator $\F_{-\frac{1}{3}}$ is again deduced by (\ref{317}):
$\F_{-\frac{1}{3}}=  E_1+\frac{1}{2^{1/3}}E_2$, and the vectors of the set $\I^{(-1/3)}$
are:  $\varphi_1^{(-1/3)}=\frac{1}{2}(1+\frac{1}{2^{1/3}},-1+\frac{1}{2^{1/3}})$,
$\varphi_2^{(-1/3)}=\frac{1}{2}(-1+\frac{1}{2^{1/3}},1+\frac{1}{2^{1/3}})$ and
$\varphi_3^{(-1/3)}=\frac{1}{2^{5/6}}(1,1)$. The dual vectors are the ones corresponding to
$-1-\alpha=-1+\frac{1}{3}=-\frac{2}{3}$, that is 
$\varphi_1^{(-2/3)}=\frac{1}{2}(1+\frac{1}{2^{2/3}},-1+\frac{1}{2^{2/3}})$,
$\varphi_2^{(-2/3)}=\frac{1}{2}(-1+\frac{1}{2^{2/3}},1+\frac{1}{2^{2/3}})$ and
$\varphi_3^{(-2/3)}=\frac{1}{2^{7/6}}(1,1)$. Again, it is straightforward to verify that the
expansions in (\ref{315}) hold true, that is that, for any $f\in \Hil$, 
$f=\sum_{i=1}^3<\varphi_i^{(-1/3)},f> \varphi_i^{(-2/3)} = \sum_{i=1}^3<\varphi_i^{(-2/3)},f>
\varphi_i^{(-1/3)}$. Moreover, we can also check successfully the estimates given in
Proposition 1: $\|f\|^2 \leq
\sum_{i=1}^3|<\varphi_i^{(-1/3)},f>|^2 \leq 2^{1/3}\|f\|^2$ and 
 $\frac{1}{2^{1/3}}\|f\|^2 \leq \sum_{i=1}^3|<\varphi_i^{(-2/3)},f>|^2 \leq \|f\|^2$.

\vspace{4mm}

\noindent
{\bf Example 2} 

We consider now another finite dimensional example a little bit more complicated. Let us
take $\Hil = \R^3$, and let us consider the four (linearly dependent) vectors
$\varphi_1=(1,0,0)$, $\varphi_2=(0,1,0)$, $\varphi_3=(0,0,1)$ and
$\varphi_4=\frac{1}{\sqrt{3}}(1,1,1)$. These vectors are a $(1,2)$-frame: $$ \|f\|^2 \leq
\sum_{i=1}^3|<\varphi_i,f>|^2 \leq 2\|f\|^2. $$ (Again, the frame bounds may be optimized).
For any $f=(f_x,f_y,f_z)\in \Hil$ we get now $\F_1f=\frac{1}{3}(4f_x+f_y+f_z,f_x+4f_y+f_z,
f_x+f_y+4f_z)$, so that  $$
\F_1\equiv \!\! \frac{1}{3}\left(
\begin{array}{ccc}
4 & 1 & 1 \\ 
1 & 4 & 1  \\ 
1 & 1 & 4  \\ 
\end{array}
\right).
$$
Using the spectral decomposition of $\F_1$ we can write
\be
\F_1=E_1+E_2+2E_3,
\label{318}
\en
where $\lambda_1=\lambda_2=1$ and $\lambda_3=2$ are the eigenvalues of $\F_1$,
$\eta_1=\frac{1}{\sqrt{2}}(-1,0,1)$,  $\eta_2=\frac{1}{\sqrt{6}}(1,-2,1)$ and
$\eta_3=\frac{1}{\sqrt{3}}(1,1,1)$ are the corresponding orthonormal eigenvectors and $E_k$
are the projectors on these vectors. We consider here only the case $\alpha=-\frac{1}{2}$.

In this case $\F_{-\frac{1}{2}} = E_1+E_2+\frac{1}{\sqrt{2}}E_3$, and the $(1,1)$-frame is
obtained as usual: $\varphi_i^{(-1/2)}= \F_{-\frac{1}{2}}\varphi_i$, for $i=1,2,3,4$.
We get $\varphi_1^{(-1/2)}= \frac{1}{3}(2+\frac{1}{\sqrt{2}},-1+\frac{1}{\sqrt{2}},
-1+\frac{1}{\sqrt{2}})$, $\varphi_2^{(-1/2)}=
\frac{1}{3}(-1+\frac{1}{\sqrt{2}},2+\frac{1}{\sqrt{2}}, -1+\frac{1}{\sqrt{2}})$, 
$\varphi_3^{(-1/2)}= \frac{1}{3}(-1+\frac{1}{\sqrt{2}},-1+\frac{1}{\sqrt{2}},
2+\frac{1}{\sqrt{2}})$ and $\varphi_4^{(-1/2)}=
\frac{1}{\sqrt{6}}(1,1, 1)$. It is easy to verify that, for any $f\in \Hil$, $f=
\sum_{i=1}^4<\varphi_i^{(-1/2)},f> \varphi_i^{(-1/2)}$ and that
$\sum_{i=1}^4|<\varphi_i^{(-1/2)},f>|^2 = \|f\|^2$.

Again, we observe that the various vectors are not normalized. Of course any attempt of
normalizing these vectors would imply the loosing of their nature of (1,1)-tight-frame.

\vspace{4mm}

\noindent
{\bf Example 3} 

We now consider an infinite dimensional example, previously discussed in \cite{dau2}. Let us
consider two real numbers $p_0, q_0\in [\pi, 2\pi[$. Let us then introduce $\lambda =
\frac{2\pi}{p_0-q_0}$ and  a $C^\infty$-function $\nu(x)$ which is zero if $x\leq 0$ and $1$
if $x\geq 1$. In \cite{dau2} it is discussed that the function
\beano
g(x) = q_0^{-1/2}\left\{
\begin{array}{ll}
0 \hspace{6cm}x\leq -\frac{\pi}{p_0}, x\geq \frac{\pi}{p_0} \\
 \sin\left[\frac{\pi}{2}\nu\left(\frac{1}{\lambda}(\frac{\pi}{p_0}+x)\right)\right]
\hspace{2cm} -\frac{\pi}{p_0} \leq x \leq -\frac{\pi}{p_0}+\lambda \\
1 \hspace{6cm} -\frac{\pi}{p_0}+\lambda \leq x \leq \frac{\pi}{p_0}-\lambda \\
 \cos\left[\frac{\pi}{2}\nu\left(\frac{1}{\lambda}(x-\frac{\pi}{p_0}+\lambda)\right)\right]
\hspace{1.5cm} \frac{\pi}{p_0}-\lambda \leq x \leq \frac{\pi}{p_0}, \\
\end{array}
\right.
\enano

satisfies the following relations: 

(a) \hspace{2cm} $|\mbox{supp } g| =\frac{2\pi}{p_0}$;

(b) \hspace{2cm} $\sum_{k\in \Z} |g(x-kq_0)|^2 =\frac{1}{q_0}.$

Let now $W(p,q)$ be the Weyl-Heisenberg operator defined by 
$$[W(p,q)f](x) \equiv e^{ipx}
f(x-q),$$
 for any square integrable function $f(x)$. Given a function $g(x) \in \Hil={\cal
L}^2(\R)$ in reference \cite{dau2} it is shown that, by means of the above properties of the
function $g(x)$, the set $\{g_{mn}\}= \{[W(mp_0,nq_0)g](x)\}$, for $m$ and $n$ both
belonging to $\Z$, is a tight frame: $\sum_{m,n\in \Z} |<g_{mn},f>|^2 =
\frac{2\pi}{p_0q_0}\|f\|^2$ for all $f\in \Hil$. In this case it is particularly simple to
write down the operator $\F_1$, which is simply $\F_1=\frac{2\pi}{p_0q_0}\Id$. Therefore,
as it is obvious, the $(1,1)$-frame is formed by the vectors $g_{mn}^{(-1/2)}=
\F_{-\frac{1}{2}}g_{mn}= \sqrt{\frac{p_0q_0}{2\pi}} g_{mn}$.

Other non-tight frames can be obtained considering different values of $\alpha$.

\vspace{5mm}

A brief remark is in order: as it is clear from the above examples our procedure, which is
quite general, can be easily implemented whenever we deal with finite dimensional Hilbert
spaces or, when the original frame is tight. In both these situations, in fact, the operator
$\F_1$ can be written in an easily handled form. On the contrary, for general frames in an
infinite dimensional space, the situation, thought being clear from the theoretical point of
view, is not as much easy to be implemented: the spectral projections must be replaced by
the spectral measure and the finite sum which defines  $\F_1$ with an
integral over a finite interval. $\F_\alpha$ is defined easily but, in general, it is not
so clear which should be the explicit expressions of the vectors of the new frames,
$\varphi_i^{(\alpha)}=\F_\alpha \varphi_i$.

There is a particular situation, which however is rather rare, in which the computations are
easily done even in the infinite dimensional situation, namely when $\F_1$ is a finite-rank
operator. In this case, in fact, the operator $\F_1$ can still be written as
$\F_1= \sum_{k=1}^N\lambda_kE_k$, where $\lambda_k$ is the k-th eigenvalue of $\F_1$ and
$E_k$ the related projection operator. Of course, it would be interesting to obtain a
condition on the original frame which ensures that the related operator $\F_1$ is of
finite-rank. For instance it is clear that to any tight frame, and in particular to any
orthonormal basis, it cannot correspond a finite-rank operator $\F_1$ since, being
in this case $\F_1$ a multiple of the identity in $\Hil$, its rank is necessarily the
whole $\Hil$.

In any case, this kind of problem is not new. It is just the same problem which we have to
face with when we try to obtain $\F_{-1}=(F^*F)^{-1}$ using the standard approach,
\cite{dau1}: except that in very few situations  one has to proceed perturbatively,
following, for instance, the procedure discussed in the Section 2.

We will now show how this perturbation approach can be used for the computation of
$\F_\alpha$, at least when the frame bounds $A$ and $B$ do not differ too much. In the next
Section we will discuss a different perturbative approach which relaxes this condition.

The problem we want to discuss is the way in which the dual operators $\F_\alpha$ and
$\F_{-1-\alpha}$ can be computed perturbatively. We will restrict to some
choices of $\alpha$ which have a particular interest. 

Let $\alpha \in \N$. Using (\ref{28}) we have $\F_\alpha=\left(\frac{A+B}{2}\right)^\alpha
(\Id-R)^\alpha$, which is already given as a (finite) sum of positive powers of the operator
$R=\Id -\frac{2}{A+B}F^*F$. More involved is the computation of $\F_{-1-\alpha}$. Introducing
the natural number $M_\alpha =1+\alpha$, we see that $\F_{-1-\alpha}=(\F_{-1})^{M_\alpha}$.
Therefore, the problem of getting an expansion for  $\F_{-1-\alpha}$, is reduced to the
analogous problem for  $\F_{-1}$, which has been already discussed in Section 2.

The situation is specular if we take $\alpha$ such that $-\alpha \in \N$: in this case 
$\F_{-1-\alpha}$ can be found exactly, while the computation of $\F_\alpha$ is again reduced
to the computation of $\F_{-1}$.

Moreover, when $\alpha=0$, we go back to the standard procedure. 

For what concerns other values of $\alpha$, we discuss here only the case
of particular interest $\alpha=-\frac{1}{2}$, which, as we have seen, gives rise to a
(1,1)-tight-frame.

Since $\|R\|\leq 1$, we can expand the operator $\F_{-\frac{1}{2}}$ as
$$
\F_{-\frac{1}{2}} = \sqrt{\frac{2}{A+B}}\sum_{k=0}^\infty 
\left(
\begin{array}{c}
-1/2  \\ 
k \\ 
\end{array}
\right)
(-R)^k.
$$
Therefore we have, for any $i\in J$,
$$
\varphi_i^{(-\frac{1}{2})} = \F_{-\frac{1}{2}} \varphi_i = \sqrt{\frac{2}{A+B}}(\Id -R)^{-1/2}
\varphi_i =  \sqrt{\frac{2}{A+B}}\sum_{k=0}^\infty
\left(
\begin{array}{c}
-1/2  \\ 
k \\ 
\end{array}
\right)
(-R)^k  \varphi_i.
$$
We define now an approximation of $\varphi_i^{(-\frac{1}{2})}$ of {\em degree} $N$:
$$
\varphi_i^{(-\frac{1}{2}),N} \equiv \sqrt{\frac{2}{A+B}}\sum_{k=0}^N
\left(
\begin{array}{c}
-1/2  \\ 
k \\ 
\end{array}
\right)
(-R)^k  \varphi_i.
$$
Given $f\in \Hil$, how large is the error that we do when we substitute
$\varphi_i^{(-\frac{1}{2})}$ with $\varphi_i^{(-\frac{1}{2}),N}$ in the exact expansion $f
= \sum_{i\in J}<\varphi_i^{(-\frac{1}{2})}, f> \varphi_i^{(-\frac{1}{2})}$? In other words,
if we define the vector $$ f^{(N)} \equiv \sum_{i\in
J}<\varphi_i^{(-\frac{1}{2}),N}, f> \varphi_i^{(-\frac{1}{2}),N},
$$
how big is the norm $\|f-f^{(N)}\|$? This norm can be estimated telescoping it twice. Let us
start introducing the self-adjoint operator
$$
T_N \equiv (\Id-R)^{-1/2}-\sum_{k=0}^N
\left(
\begin{array}{c}
-1/2  \\ 
k \\ 
\end{array}
\right)
(-R)^k. 
$$
Therefore, $\sqrt{\frac{2}{A+B}} T_N \varphi_i = \varphi_i^{(-\frac{1}{2})} -
\varphi_i^{(-\frac{1}{2}),N}$, for all $i\in J$. It is easy to see that
\beano
&&\|f-f^{(N)}\| \leq \|\sum_{i\in J}
<(\varphi_i^{(-\frac{1}{2})}-\varphi_i^{(-\frac{1}{2}),N}), f> \varphi_i^{(-\frac{1}{2})}\| +
\\
 &&+\|\sum_{i\in J}
<(\varphi_i^{(-\frac{1}{2}),N}-\varphi_i^{(-\frac{1}{2})}), f>
(\varphi_i^{(-\frac{1}{2})}-\varphi_i^{(-\frac{1}{2}),N})\| + \\
&&+ \|\sum_{i\in J}
<\varphi_i^{(-\frac{1}{2})}, f>
(\varphi_i^{(-\frac{1}{2})}-\varphi_i^{(-\frac{1}{2}),N})\|,
\enano
which shows that the operator $T_N$ plays a relevant role in the estimate of this norm. Let
us consider in some details, for instance, the computation of the first contribution above. 
We first notice that $<(\varphi_i^{(-\frac{1}{2})}-\varphi_i^{(-\frac{1}{2}),N}), f> =
\sqrt{\frac{2}{A+B}}<T_N\varphi_i, f>= \sqrt{\frac{2}{A+B}}<\varphi_i,T_N f>$. Moreover,
using the notation of Section 2, we also deduce that
$<g,\varphi_i^{(-\frac{1}{2})}>=<g,\F_{1/2}\tilde \varphi_i>=<\F_{1/2}
g, \tilde \varphi_i>$.  Therefore, this contribution can be estimated in the
following fashion: 
\beano 
\hspace{-1cm} &&\| \sum_{i\in J}
<(\varphi_i^{(-\frac{1}{2})}-\varphi_i^{(-\frac{1}{2}),N}), f> \varphi_i^{(-\frac{1}{2})}\|
=  \sup_{\|g\|\leq 1} \left| \sqrt{\frac{2}{A+B}} \sum_{i\in J}
<\varphi_i, T_N f> <\F_{1/2}g, \tilde \varphi_i>\right|= \\
&& = \sqrt{\frac{2}{A+B}}\sup_{\|g\|\leq 1}  \left|<\F_{1/2}g, T_N f>\right| = 
\sqrt{\frac{2}{A+B}} \|\F_{1/2}T_N f\| \leq \sqrt{\frac{2B}{A+B}} \|T_N\|\,
\|f\|.
\enano
Here we have used, among the other things, the fact that: $\|F_{1/2}\|\leq
\sqrt{B}$, which directly follows from (\ref{35}). The other contributions above can be
estimated in analogous way. The final result is the following:
$$
\|f-f^{(N)}\| \leq \left(2\sqrt{\frac{2B}{A+B}}+\frac{2B}{A+B}\|T_N\|\right)\|T_N\|\, \|f\|.
$$
We see here the necessity of estimating $\|T_N\|$. The details of this computation are given
in the Appendix, and the result is in formula (\ref{A2}). Here we give only the result:
$$
\|T_N\| \leq \left(\frac{B-A}{2A}\right)^{N+1} \sqrt{\frac{A+B}{2A}}.
$$
The conclusion is therefore the following inequality:
$$
\|f-f^{(N)}\| \leq \sqrt{\frac{B}{A}}\left(\frac{B-A}{2A}\right)^{N+1}\left(2+ 
\sqrt{\frac{B}{A}} \left(\frac{B-A}{2A}\right)^{N+1}\right),
$$
which converges to zero in the same fashion of the analogous estimates given in \cite{dau2},
at least if $B<3A$. Obviously the speed of this convergence increases as much as $A$
becomes closer to $B$.

\vspace{5mm}

We want to end this long Section with some {\em reversibility} remarks. The first question
is the following: given a frame, it gives rise to a single operator $\F_1=F^*F$. Is this
correspondence a bijection? That is, given a self-adjoint operator $\F_1$ which satisfies the
properties discussed at the beginning of this Section, is it possible to associate to this
operator a single frame? 

It is easy to find counterexamples which show that in general this is not so. The first
trivial case is the following: if the set $\{\varphi_i\}$ is a frame then also
$\{e^{i \alpha_i} \varphi_i\}$ is a frame for any choice of the real constants $\alpha_i$,
and both the frames correspond to the same $\F_1$ operator. Therefore it is not possible,
given $\F_1$, to fix a single 'generating frame'. However, it could happen, in principle,
that $\F_1$ fixes only a frame of {\em rays} in the Hilbert space, and not a frame of
vectors. This cannot happen, in general, for finite dimensional Hilbert spaces as we can
show with a very easy example: let us suppose that dim$(\Hil)=2$ and that $\Hil$ is a real
Hilbert space. Then $\F_1$ is a $2\times 2$ real and symmetric matrix with three independent
entries. Let us now consider a frame of $\Hil$ of four vectors. Condition
$\F_1f=\sum_{i=1}^4<\varphi_i,f>\varphi_i$ for any $f\in \Hil$ gives three independent
equations which become seven if we also consider the normalization conditions of the
vectors. They are not enough to fix all the vectors of the frame! 

The question is still open for infinite dimensional Hilbert spaces.

The second remark is very much related to the previous one: given an $(A,B)$-frame
$\{\varphi_i\}$ and its related operator $\F_1$, we find an unique $(1,1)$-frame
$\{\varphi_i^{(-1/2)}\}$. The question is again a problem of reversibility: given now an 
$(1,1)$-frame $\{\phi_i\}$, is it possible to find an unique $(A,B)$-frame $\{\Psi_i\}$ such
that $\phi_i = \Psi_i^{(-1/2)}$? Again, the answer is in general negative and it can be
shown simply by giving examples of $(1,1)$-frames each of which can be obtained by different
non-tight frames. Before discussing these examples we remind that normalization of the
vectors of the frame is very important: it is enough to think to the examples of
$(1,1)$-frames discussed in this Section, which are certainly not orthonormal bases due to
the lack of normalization of their vectors. It is well known, in fact, that such a frame
turns out to be an orthonormal basis if its vectors are all normalized, \cite{dau1}.

We start observing the following property: 

- If $\{\varphi_i\}$ is a frame then also $\{\Phi_i= \alpha \varphi_i\}$ is a frame and we
have $\varphi_i^{(-1/2)}=\Phi_i^{(-1/2)}$ for all $i\in J$ and $\forall \alpha>0$.

The proof of this statement is a straightforward computation. With obvious notation we call
$F_\Phi$ and $F_\varphi$ the frame operators related to the two different frames. By the
definition we have, for any $f\in \Hil$:
$$
F_\Phi^*F_\Phi f=\sum_{i\in J}<\Phi_i,f>\Phi_i=\alpha^2 \sum_{i\in J}<\varphi_i,f>\varphi_i =
\alpha^2 F_\varphi^*F_\varphi f,
$$
so that $F_\Phi^*F_\Phi = \alpha^2 F_\varphi^*F_\varphi$. Therefore we have
$$
\Phi_i^{(-1/2)}= (F_\Phi^*F_\Phi)^{-1/2} \Phi_i = \frac{1}{\alpha}
(F_\varphi^*F_\varphi)^{-1/2} \alpha \varphi_i =  \varphi_i^{(-1/2)}.
$$

This is a particular case of a more general situation which is easily considered if dim
$\Hil <\infty$. Let $\{\varphi_i\}$ be a given frame, with $i=1,...,N$, and let $F_\varphi$
be its frame operator. Let us suppose that it exists a self-adjoint operator $A\in B(\Hil)$
which commutes with $F_\varphi^*F_\varphi$ and for which there exist two positive real
numbers $0<A_-\leq A_+<\infty$ such that  $A_-\Id \leq  A\leq A_+\Id$. In general
any operator of the form $\alpha \Id$ satisfies these conditions for any strictly positive
$\alpha$, independently of the form of $\F_\varphi$, and this is the case we have considered
above, as we will see in a moment. In general such an operator, if it exists, must be
related to the frame operator. Let us now define $\Phi_i
=A^{1/2}\varphi_i$ for all $i=1,...,N$. The set of these new vectors is a frame whose
related frame operator is $F_\Phi = F_\varphi A^{1/2}$. We have: $$
\Phi_i^{(-1/2)}=(F_\Phi^*F_\Phi)^{(-1/2)} \Phi_i=(F_\varphi^*F_\varphi)^{(-1/2)} A^{-1/2} 
A^{1/2}\varphi_i=(F_\varphi^*F_\varphi)^{(-1/2)} \varphi_i= \varphi_i^{(-1/2)}, $$
as it was to be shown. The extension to the case dim $\Hil =\infty$ contains some subtle
points and it will not be considered here.

\section{A Perturbation Approach}

We have discussed in the previous Sections the difficulties of finding the operator
$\F_{-1}=(F^*F)^{-1}$ for a given frame, in general. The way in which this problem is
overcame in the literature is to consider a perturbative approach whose speed of
convergence, however, depends strongly on the difference of the frame bounds, $B-A$: the
smaller this quantity the more rapid the convergence of the procedure. In this Section
we will propose a different perturbation approach which relax this condition, so that it
can be also used when $B-A$ is not so small.

The idea is very simple and it is strongly related to the fact that the operator $\F_1=F^*F$
is strictly positive. Let $\{\varphi_i\}$ be a frame and $F$ its associated frame operator.
We know that
\be
A\Id \leq F^*F \equiv \F_1 \leq B\Id.
\label{41}
\en
In all this Section we will not assume that $\frac{B-A}{B+A}$ is much smaller than $1$. On
the contrary, the procedure we are going to discuss is particularly useful in the case in
which this ratio approaches the unity, that is in the case in which $B$ is much larger than
$A$.  How can we construct a perturbative approach to get the operator $\F_{-1} =
(F^*F)^{-1}$? We will propose a procedure which depends on the values of the frame bounds.

\vspace{4mm}

Let us first assume that
\be
1<A\leq B <\infty.
\label{42}
\en
As it was discussed in the previous Section, due to the nature of the operator $\F_1$ it is
possible to find a spectral decomposition for this operator and we can build up functions of
this operator. For instance, recalling equation (\ref{31})
\be
<\F_1\Phi,\Psi> = \int_A^B\lambda \, d<E_\lambda\Phi, \Psi>, \hspace{2cm} \forall \Phi, \Psi
\in \Hil,
\label{43}
\en
we see that the operator $\Log_B(\F_1)$ is defined by
\be
<\Log_B(\F_1)\Phi,\Psi> = \int_A^B\Log_B(\lambda) \, d<E_\lambda\Phi, \Psi>, \hspace{2cm}
\forall \Phi, \Psi \in \Hil.
\label{44}
\en
The integral above is surely well behaved since $\lambda \geq A >1$. Using functional
calculus it is easy to verify that inequality (\ref{41}) gives rise to
\be
\alpha \Id \leq \Log_B(\F_1) \leq \Id.
\label{45}
\en
This result is, of course, a consequence of the fact that $A\Id$ and $B\Id$ both commute
with $\F_1$. In (\ref{45}) we have defined $\alpha = \Log_B(A)$. Of course $0<\alpha\leq
1$. The relevant difference with respect to the original double inequality (\ref{41}) is
that the two new bounds $\alpha$ and $1$, are necessarily near to each other. We define, in
analogy with the standard approach, an operator $R$ via
\be
R=\Id -\frac{2}{1+\alpha}\Log_B(\F_1).
\label{46}
\en
Using the bounds (\ref{45}) we find that the operator $R$ satisfies the following
double inequality:
\be
-\frac{1-\alpha}{1+\alpha} \Id \leq R \leq \frac{1-\alpha}{1+\alpha}\Id.
\label{47}
\en
Equation (\ref{46}) can now be inverted and it gives, after some minor calculation,
\be
\F_{-1} = (F^*F)^{-1} = \frac{1}{\sqrt{A\, B}} B^{\frac{1+\alpha}{2}R} = \frac{1}{\sqrt{A\,
B}} e^{\ln B \frac{1+\alpha}{2} \, R}.
\label{48}
\en
Let us now see how the perturbative approach works. We call, as usual, $\tilde \varphi_j =
\F_{-1} \varphi_j$, which is the {\em exact} dual frame, and we define its N-th
approximation by
\be
\tilde \varphi_j^N \equiv \frac{1}{\sqrt{A\, B}} \sum_{k=0}^N \frac{1}{k!}
\left(\ln B \frac{1+\alpha}{2} R \right)^k \, \varphi_j, \hspace{2cm} \forall j\in J.
\label{49}
\en
We are interested in computing which error do we do if, instead of $\sum_{i\in J}
<\varphi_i, f> \tilde \varphi_i = f$, we compute $f^{(N)} \equiv \sum_{i\in J} <\varphi_i,
f> \tilde \varphi_i^N$. In other words, we want to estimate $\|f-f^{(N)}\|$. The procedure is
similar to the one we have discussed in the previous Section. Let us define the
self-adjoint operator $Z_N$, $$
Z_N \equiv e^{\ln B \frac{1+\alpha}{2} \, R} - \sum_{k=0}^N \frac{1}{k!}
\left( \ln B \frac{1+\alpha}{2} \, R \right) ^k.
$$
Using the definitions and the following equality
$$
<\tilde \varphi_j- \tilde \varphi_j^N,g> =\frac{1}{\sqrt{AB}}<Z_N\F_1\tilde \varphi_j,g> = 
\frac{1}{\sqrt{AB}}<\tilde \varphi_j,\F_1 Z_N g>
$$
where, as usual, $\tilde \varphi_j = \F_{-1} \varphi_j$, we have
\beano
&&\|f-f^{(N)}\| = \sup_{\|g\|\leq 1} \left| \sum_{j\in J} <f,\varphi_j> <\tilde \varphi_j- \tilde
\varphi_j^N,g>\right|= \\
&&= \frac{1}{\sqrt{AB}}\sup_{\|g\|\leq 1} \left| <f,\F_1 Z_N g> \right| \leq
\sqrt{\frac{B}{A}} \|f\| \, \|Z_N\|.
\enano
The problem is, therefore, to evaluate $\|Z_N\|$. In the Appendix we will deduce that 
$$
\|Z_N\| \leq \sqrt{\frac{B}{A}}\frac{1}{(N+1)!}\left( \frac{1-\alpha}{2} \ln B\right) ^{N+1},
$$
so that 
\be
\|f-f^{(N)}\| \leq \frac{B}{A} \, \|f\| \, \frac{\left( \frac{1-\alpha}{2} \ln B\right)
^{N+1}}{(N+1)!}, \label{410}
\en
which shows that $\|f-f^{(N)}\|$ is exponentially convergent to zero when  $N\rightarrow
\infty$, independently of the values of $A$ and $B$, which can also be very different from
each other.

\vspace{5mm}

Let us now assume that 
\be
0<A\leq B <1.
\label{411}
\en
In this case the inequality (\ref{41}) can be rewritten as 
\be
\beta \Id \leq \Log_A(\F_1) \leq  \Id,
\label{412}
\en
where $\beta \equiv \Log_A B$ is strictly positive and $\beta \leq 1$. We define the
operator $R$ in a way which is similar to what we have done in (\ref{46}):
\be
R=\Id -\frac{2}{1-\beta}\Log_A(\F_1).
\label{413}
\en
The operator $R$ is bounded, $\|R\| \leq \frac{1-\beta}{1+\beta}$, and the operator $\F_{-1}$
differs from the one given in (\ref{48}) and turns out to be 
\be
\F_{-1} =  \frac{1}{\sqrt{A\, B}} A^{\frac{1+\beta}{2}R} = \frac{1}{\sqrt{A\,
B}} e^{\ln A\frac{1+\beta}{2} \, R}.
\label{414}
\en
Again, if we define $f^{(N)}$ in the usual way, we get for $\|f-f^{(N)}\|$ an analogous 
estimate as the one we have obtained before, showing its exponential convergence to zero
for $N$ diverging.

\vspace{4mm}

The last condition on the frame bounds $A$ and $B$ is
\be
0< A \leq 1, \hspace{3cm} B\geq 1,
\label{415}
\en
which can be reconduced to the first situation simply by multiplying (\ref{415}) for, say,
$\frac{2}{A}$. In this way this inequality can be rewritten as 
\be
2\Id \leq \frac{2}{A} \F_1 \leq b\Id,
\label{416}
\en
where $b\equiv \frac{2B}{A}$. At this point we call $\delta \equiv \Log_b 2$, so that
(\ref{416}) becomes
\be
\delta \Id \leq \Log_b(\frac{2}{A} \F_1) \leq \Id,
\label{417}
\en
and we can proceed exactly as we have done before. We obtain, after some minor computation,
\be
\F_{-1} = \frac{1}{\sqrt{A\, B}} b^{\frac{1+\delta}{2} R} = \frac{1}{\sqrt{A\, B}} 
e^{\ln b \frac{1+\delta}{2} R }.
\label{418}
\en
The estimate of $\|f-f^{(N)}\|$ goes through as usual. Again, this quantity converges to zero
exponentially with $N$.

The following remark is in order: while in the standard approach the zero
approximation of $\tilde \varphi_j$, $\tilde \varphi_j^0$, is nothing but $\frac{2}{A+B} 
\tilde \varphi_j$, with this approach the zero-th approximation turns out to be
$\frac{1}{\sqrt{A\,B}}  \tilde \varphi_j$. The arithmetic mean is replaced with the geometric
one.

Of course the standard perturbative approach we have used in the previous Section in order to
discuss how to obtain $\F_\alpha$ could be replaced with the different approach we have
just developed. However, we will not consider this analysis here.

\section{Frames Versus Orthonormal Bases: Some Remarks}

This Section, which concludes the paper, is devoted to discuss some differences
between a set of orthonormal vectors and a frame. 

As we know, \cite{ric}, a set $\{\varphi_j\}$, $j\in J$, is an orthonormal basis of a
given Hilbert space $\Hil$, if and only if one of the following equivalent conditions is
satisfied:

(i) for all $f\in \Hil$ then $\|f\|^2 = \sum_{j\in J} |<f, \varphi_j>|^2$;

(ii) for all $f\in \Hil$ then $f = \sum_{j\in J} <\varphi_j, f> \varphi_j$;

(iii) for all $f, g\in \Hil$ then $<f, g> = \sum_{j\in J}  <f, \varphi_j> <\varphi_j, g>$;

(iv) if $f\in \Hil$ is such that $<f,\varphi_j>=0$ for all $j\in J$, then $f=0$.

Similar conditions can also be restated for an $(A,B)$-frame $\{\phi_j\}$ with frame
operator $F$. Introducing the dual frame of vectors $\tilde \phi_j = (F^*F)^{-1} \phi_j$,
the above conditions become:

(i$'$) for all $f\in \Hil$ then $A \|f\|^2 \leq  \sum_{j\in J} |<f, \phi_j>|^2 \leq B 
\|f\|^2$;

(ii$'$) for all $f\in \Hil$ then $f = \sum_{j\in J} <\phi_j, f> \tilde \phi_j$;

(iii$'$) for all $f, g\in \Hil$ then $<f, g> = \sum_{j\in J}  <f, \phi_j> <\tilde \phi_j,
g>$;

(iv$'$) if $f\in \Hil$ is such that $<f,\phi_j>=0$ for all $j\in J$, then $f=0$.

These conditions can also be rewritten in terms of the frame operator $F$. In this way
conditions ii$'$) and iii$'$) collapse and the conditions are:

(a) $A\Id \leq F^*F \leq B\Id$;

(b)  the operator $(F^* F)$ can be inverted and its inverse $(F^* F)^{-1}$ is bounded in
$\Hil$;

(c) ker $F = \{ 0\}$.

The equivalence between (a) and (i$'$) is well established in the literature. For what
concerns the equivalence of condition (b) with the conditions (ii$'$) and (iii$'$), this is
easy to be deduced: since the operator $(F^* F)^{-1}$ belongs to $B(\Hil)$, we can
 define a new operator $\tilde F = F (F^* F)^{-1}$ which maps $\Hil$ into ${\Bbb
l}^2(J)$ such that $\tilde F^* F = F^* \tilde F=\Id$. Defining $\tilde \phi_j \equiv (F^*
F)^{-1}\phi_j$ one direction of the equivalence is proved. The converse implication is
immediate, since the existence of the dual frame already implies the existence of
$(F^*F)^{-1}$ and its boundedness. Finally, if ker $F=\{0\}$, then the operator $F$ is
injective, so that $Ff=0$ only if $f=0$. Vice-versa, if $Ff=0$ implies that $f=0$, then
$F$ is injective as an operator from $\Hil$ into ${\Bbb l}^2(J)$, so that its kernel
contains only the zero vector.

It is
widely discussed in the literature that condition (a) implies condition (b), or, in other
terms, that (i$'$) implies (ii$'$) and (iii$'$), essentially as a consequence of the
inequality $F^*F\geq A\Id$. It is also very easy to see that condition (i$'$) implies
condition (iv$'$) or that (a) implies (c): let us suppose that $f\in \Hil$ is orthogonal to
all the frame vectors. If condition (i$'$) holds it is obvious that $\|f\|=0$, so that
$f=0$.

The main point is now to show that, contrarily to what happens for orthonormal bases, the
converse implications are not all true. In particular, we will show that condition (b)
implies condition (a), while condition (c) is not equivalent to the other two conditions.
In order to show the implication $(b) \Rightarrow (a)$ we observe first that the operator 
$(F^* F)^{-1}$ is positive as its inverse. Therefore, using the fact that if the bounded
operator $T$ is positive then $T\|T\| \geq T^2$, \cite{brat}, we deduce that the operator
$F^*F$ satisfies the following inequalities: $\frac{1}{\|(F^*F)^{-1}\|}\Id \leq F^*F \leq
\|F^*F\|\Id.$ Therefore $F$ is a frame operator.

The fact that the last condition does not imply, say, the first one, follows from the fact
that the condition of $F$ being injective cannot of course imply also its surjectivity.
Therefore the set of vectors related to this $F$ in general does not generate the whole
Hilbert space, so that it is not a frame, in general. On the other hand, if we also require
to $F$ of being surjective, then it would be invertible, condition which we know is not
required at all to any frame operator. This is what happens for Riesz bases which, as we
have discussed in the Introduction, generate the whole Hilbert space  and contain
only linearly independent vectors, \cite{cdf}. 

We conclude, therefore, noticing that, even if the conditions which define an orthonormal
basis can be easily extended to general frames, these new conditions are not in general all
equivalent. This difference is essentially due to the linear dependence of the
vectors of the frame.

\vspace{1cm}

\noindent{\large \bf Acknowledgments} \vspace{5mm}

This work has been supported by M.U.R.S.T.

\vspace{2cm}

 \appendix

\renewcommand{\theequation}{\Alph{section}.\arabic{equation}}

 \section{\hspace{-.7cm}ppendix :  Norm Estimates}

This Appendix is devoted to show the details of the estimate of the norms of a couple of
operators which are used in the paper. 

Let us start with the first operator,
$$
T_N \equiv (\Id-R)^{-1/2}-\sum_{k=0}^N
\left(
\begin{array}{c}
-1/2  \\ 
k \\ 
\end{array}
\right)
(-R)^k,
$$
which appears in Section 3 in the computation of $\|f-f^{(N)}\|$. The estimate proceeds
using the functional calculus in the following way: let us introduce the spectral family
$E_\lambda$ of the bounded and self-adjoint operator  $R$. Of course this spectral family
coincides with the one of the operator $\F_1=F^*F$, see equality (\ref{28}). As in 
equation (\ref{32}), we can write \be
<R \Phi,\Psi> = \int_{-\gamma}^{\gamma}\lambda \: d<E_\lambda\Phi, \Psi>,
\hspace{2cm} \forall \Phi, \Psi \in \Hil,
\label{A1}
\en
where $\gamma \equiv \frac{B-A}{B+A}$. Being $(\Id -R)$ a positive operator bounded from below
from a strictly positive constant, we can consider the self adjoint operator $(\Id -R)^{-1/2}$,
\cite{ric}, which can be written, in terms of the spectral family $E_\lambda$, as 
$$
<(\Id -R)^{-1/2} \Phi,\Psi> = \int_{-\gamma}^{\gamma}(1 -\lambda)^{-1/2} \,
d<E_\lambda\Phi, \Psi>, \hspace{2cm} \forall \Phi, \Psi \in \Hil.
$$
Obviously we have also 
$$
<\sum_{k=0}^N
\left(
\begin{array}{c}
-1/2  \\ 
k \\ 
\end{array}
\right)
(-R)^k \Phi,\Psi> = 
\int_{-\gamma}^{\gamma}\sum_{k=0}^N
\left(
\begin{array}{c}
-1/2  \\ 
k \\ 
\end{array}
\right)
(-\lambda)^k \,
d<E_\lambda \Phi,\Psi>.
$$
Therefore, for any $f\in \Hil$, we get
$$
\|T_Nf\|^2= \int_{-\gamma}^{\gamma}\left| (1 -\lambda)^{-1/2}- \sum_{k=0}^N
\left(
\begin{array}{c}
-1/2  \\ 
k \\ 
\end{array}
\right)
(-\lambda)^k \right|^2 \, d<E_\lambda f,f>.
$$
The problem is now reduced to estimate the following function:
$$
\rho (\lambda) = (1 -\lambda)^{-1/2}- \sum_{k=0}^N
\left(
\begin{array}{c}
-1/2  \\ 
k \\ 
\end{array}
\right)
(-\lambda)^k.
$$
The idea is to observe that $\rho (\lambda)$ is nothing but the N-th (Lagrange) remainder of the
Maclaurin approximation of the function $\Psi(\lambda)=(1 -\lambda)^{-1/2}$, so that we can
estimate $\rho (\lambda)$ in the usual way: $\rho (\lambda) = \frac{\Psi^{N+1}(\xi)}{(N+1)!}
\lambda^{N+1}$. Here $\xi$ is a suitable real belonging to $[0,\lambda]$ if $\lambda >0$ or
to  $[\lambda,0]$ if $\lambda <0$. 

With some easy estimates we see that $\Psi^{(N)}(\xi) \leq N! \left(
\frac{1}{1-\gamma}\right) ^{\frac{2N+1}{2}}$, which implies that $|\rho (\lambda)| \leq
\left( \frac{1}{1-\gamma} \right)^\frac{2N+3}{2} \gamma^{(N+1)}$. Therefore we get
$$
\|T_Nf\|^2 \leq \left( \frac{1}{1-\gamma} \right)^{2N+3} \gamma^{2(N+1)} \int_{-\gamma}^{\gamma}
\, d<E_\lambda f,f> = \left( \frac{1}{1-\gamma} \right)^{2N+3} \gamma^{2(N+1)} \|f\|^2,
$$
and, finally, 
\be
\|T_N\| \leq \frac{1}{\sqrt{1-\gamma}}\left( \frac{\gamma}{1-\gamma} \right)^{N+1}=
\left(\frac{B-A}{2A}\right)^{N+1} \sqrt{\frac{B+A}{2A}}.
\label{A2}
\en

\vspace{5mm}

The second quantity we want to estimate in this Appendix is the norm of the operator $Z_N$
introduced in Section 4,
$$
Z_N \equiv e^{\ln B\frac{1+\alpha}{2} \, R} - \sum_{k=0}^N \frac{1}{k!}
\left( \ln B\frac{1+\alpha}{2} \, R \right) ^k.
$$
The approach is exactly the same as for the previous norm: we use functional calculus,
Maclaurin expansion and the Lagrange formula for the remainder. We do not give all the
details here, but only the result which is the following: 
\be
\|Z_N\| \leq \frac{1}{(N+1)!}e^{\ln B \, \frac{1+\alpha}{2} \|R\|} \left(
\frac{1+\alpha}{2} \|R\| \, \ln B \right) ^{N+1} \leq \frac{1}{(N+1)!}
\sqrt{\frac{B}{A}}\left( \frac{1-\alpha}{2} \ln B \right)^{N+1}, 
\label{A3}
\en
where $\alpha =\Log_B A$ and, due to (\ref{47}), $\|R\| \leq \frac{1-\alpha}{1+\alpha}$.

\end{document}